\documentclass[10pt,onecolumn]{article}
\usepackage{multicol}
\usepackage[greek,english]{babel}
\usepackage[vcentering,dvips]{geometry}
\usepackage{amsfonts}
\usepackage{amssymb}
\usepackage{eucal}
\usepackage{amsxtra}
\usepackage{color}
\usepackage{setspace}
\usepackage{times}
\usepackage{titlesec}
\usepackage{sectsty}
\usepackage{fancyhdr}
\usepackage{fancyheadings}
\allsectionsfont{\large}
\usepackage{url}
\usepackage{doi}
\usepackage{hyperref}
\hypersetup{
    colorlinks,%
    citecolor=black,%
    filecolor=black,%
    linkcolor=black,%
    urlcolor=black
}

\addtocounter{section}{0} \numberwithin{equation}{section}
\newtheorem{theorem}{Theorem}[section]
\newtheorem{proposition}[theorem]{Proposition}
\newtheorem{lemma}[theorem]{Lemma}

\newcommand{\qed}{{$\hfill \Box$}}

\onehalfspace

\hypersetup{urlcolor=blue }

\begin{document}
\pagestyle{fancy}
\fancyhead[LE,RO]{\thepage}
\renewcommand{\sectionmark}[1]{\markright{\emph{ \thesection.\ #1}}{}}
\fancyhead[RE]{\leftmark}
\fancyhead[LO]{\rightmark}
\cfoot{}
\title{\textbf{The structure Jacobi operator and the shape operator of real hypersurfaces in $\mathbb{C}P^{2}$ and$\mathbb{C}H^{2}$}}

\author{\textsc{Konstantina Panagiotidou}\\
\textsc{Mathematics Division-School of Technology,}\\ \textsc{Aristotle University of Thessaloniki, Thessaloniki 54124, Greece}\\
e-mail:kapanagi@gen.auth.gr}

\date{}
\maketitle

\begin{flushleft}
\small{\textsc{Abstract:}
The aim of the paper is to present two results concerning real hypersurfaces in $\mathbb{C}P^{2}$ and $\mathbb{C}H^{2}$. More precisely, it is proved that real hypersurfaces equipped with structure Jacobi operator satisfying condition $\mathcal{L}_{X}l=\nabla_{X}l$, where \emph{X} is a vector field orthogonal to structure vector field $\xi$, do not exist. Additional real hypersurfaces equipped with shape operator \emph{A} satisfying relation $\mathcal{L}_{X}A=\nabla_{X}A$, where \emph{X} is a vector field orthogonal to $\xi$, do not exist.}
\end{flushleft}
\begin{flushleft}
\small{\emph{Keywords}: Real hypersurface, Structure Jacobi operator, Shape Operator, Lie Derivative, Complex projective space, Complex hyperbolic space.\\}
\end{flushleft}
\begin{flushleft}
\end{flushleft}

\rhead[\centering{Konstantina Panagiotidou}]{\thepage}
\lhead[\thepage]{\centering{Structure Jacobi Operator and Shape Operator of Real Hypersurfaces in $\mathbb{C}P^{2}$ and
$\mathbb{C}H^{2}$}}

\section{\textbf{Introduction}}

\medskip
A \emph{complex space form} is an n-dimensional Kaehler manifold	of constant holomorphic sectional curvature \emph{c} and it is denoted by $M_{n}(c)$. A complete and simply connected complex space form is complex analytically isometric to a complex projective space $\mathbb{C}P^{n}$, a complex Euclidean space $\mathbb{C}^{n}$ or a complex hyperbolic space $\mathbb{C}H^{n}$ if $c>0, c=0$ or $c<0$ respectively.

Let \emph{M} be a real hypersurface in non-flat complex space form $M_{n}(c)$, $c\neq0$. Then an almost contact metric structure $(\varphi, \xi, \eta, g)$ can be defined on \emph{M} induced from the Kaehler metric \emph{G} and complex structure \emph{J} on $M_{n}(c)$. The \emph{structure vector field} $\xi$ is called \emph{principal} if $A\xi=\alpha\xi$, where \emph{A} is the shape operator of \emph{M} and $\alpha=\eta(A\xi)$ is a smooth function. A real hypersurface is said to  be a \emph{Hopf hypersruface}, if $\xi$ is principal.

The study of real hypersurfaces in $M_{n}(c)$, $c\neq0$ is a classical problem in the area of Differential Geometry. In \cite{T1}, \cite{T2} Takagi  was the first who studied and classified homogeneous real hypersurfaces in $\mathbb{C}P^{n}$ and showed that they could be divided into six types, namely $A_{1}$, $A_{2}$, $B$, $C$, $D$ and $E$. In the case of $\mathbb{C}H^{n}$ Berndt (see \cite{Ber}) classified real hypersurfaces
with constant principal curvatures, when $\xi$ is principal. Such real hypersurfaces are homogeneous. Recently, Berndt and Tamaru in  \cite{BT} have given a complete classification of homogeneous real hypersurfaces in $\mathbb{C}H^{n}$, $n\geq2$.

The \emph{structure Jacobi operator} plays an important role in the study of real hypersurfaces in complex space form. It is denoted by $l$ and is given by the relation $lX=R(X,\xi)\xi$. Various results concerning different types of parallelness of $l$ have been established. Ortega, Perez and Santos in \cite{OPS} proved the non-existence of real hypersurfaces in non-flat complex space form with a parallel structure Jacobi operator, i.e. $\nabla_{X} l=0,\;X\;\in\;TM$. Perez, Santos and Suh in \cite{PSaSuh} continuing the work of \cite{OPS} considered the weaker condition $\nabla_{X}l=0$ for any vector field $X$ orthogonal to $\xi$, ($\mathbb{D}$-parallelness). They proved the non-existence of such real hypersurfaces in $\mathbb{C}P^{n}$, $n\geq3$. The condition of $\xi$-parallel structure Jacobi operator, i.e. $\nabla_{\xi}l=0$, has been studied in combination with other ones (\cite{KL}, \cite{KPSaSuh}, \cite{KK}, etc).

Conditions concerning the Lie derivative of the structure Jacobi operator is another issue that has been studied extensively. More precisely, in \cite{PSa} the non-existence of real hypersurfaces in $\mathbb{C}P^{n}$, ($n\geq3$), whose Lie derivative of the structure Jacobi operator with respect to any vector field $X$ vanishes, i.e. $\mathcal{L}_{X}l=0$, is proved. On the other hand, real hypersurfaces in $\mathbb{C}P^{n}$, $n\geq3$, whose Lie derivative of the structure Jacobi operator with respect to $\xi$ vanishes, i.e. $\mathcal{L}_{\xi}l=0$, are classified (see \cite{PSaSuh1}). Ivey and Ryan in \cite{IR} extend some of the above results in $\mathbb{C}P^{2}$ and $\mathbb{C}H^{2}$. More precisely, they proved that in $\mathbb{C}P^{2}$ and $\mathbb{C}H^{2}$ no real hypersurfaces satisfying condition $\mathcal{L}_{X}l=0$ for any vector field $X$ exists, but real hypersurfaces satisfying condition $\mathcal{L}_{\xi}l=0$ exist and they classified them. Additional, they proved that no real hypersurfaces in $\mathbb{C}P^{n}$ or $\mathbb{C}H^{n}$, $n\geq3$, satisfying condition $\mathcal{L}_{X}l=0$, for any vector field $X$ exist.

Another condition concerning the structure Jacobi operator which has been studied is $\mathcal{L}_{\xi}l=\nabla_{\xi}l$. More precisely, Perez and Santos in \cite{PSa1} classified real hypersurfaces in $\mathbb{C}P^{n}$, $n\geq3$, and the author with Ph. J. Xenos in \cite{PX1} classified real hypersrufaces in $\mathbb{C}P^{2}$ and $\mathbb{C}H^{2}$ equipped with structure Jacobi operator satisfying the latter condition.

Apart from the structure Jacobi operator, the shape operator of real hypersurface also plays an important role in the study of them. Different types of parallelness for the shape operator such as $\eta$-parallelness, pseudo-parallelness etc have been studied (\cite{LoOr}, \cite{KMa1989}, \cite{KonLoo}). In \cite{NR1} the non-existence of real hypersurfaces in complex space form with parallel shape operator, i.e. $\nabla_{X}A=0,\;X\;\in\;TM$ is referred. Additionally, Kimura and Maeda in \cite{KM}  classified real hypersurfaces in complex projective space whose shape operator is $\xi$-parallel, i.e. $\nabla_{\xi}A=0$. Finally, the Lie derivative of it has also been studied. Ki, Kim and Lee \cite{KKL1} classified real hypersurfaces in complex space forms whose shape operator is Lie $\xi$-parallel, i.e. $\mathcal{L}_{\xi}A=0$.

From the work that so far has been done, the general case of studying real hypersurfaces in complex space form $M_{n}(c)$, $c\neq0$, whose Lie derivative of a tensor field with respect to a vector field $X$ orthogonal to $\xi$ coincides with the covariant derivative of it in the same direction raised naturally.

More precisely, in the present paper real hypersurfaces in $\mathbb{C}P^{2}$ and $\mathbb{C}H^{2}$, whose structure Jacobi operator satisfies the relation
\begin{eqnarray}\label{Structure-Jacobi-Operator}
\mathcal{L}_{X}l=\nabla_{X}l,\;\mbox{where \emph{X} $\in$ $\mathbb{D}$}
\end{eqnarray}
are studied and the following theorem is proved
\medskip
\begin{theorem}\label{theorem-1}
There do not exist real hypersurfaces in $\mathbb{C}P^{2}$ or $\mathbb{C}H^{2}$, whose structure Jacobi operator satisfies relation (\ref{Structure-Jacobi-Operator}).
\end{theorem}

Additional real hypersurfaces in $\mathbb{C}P^{2}$ and $\mathbb{C}H^{2}$, whose shape operator satisfies relation
\begin{eqnarray}\label{Shape-Operator}
\mathcal{L}_{X}A=\nabla_{X}A,\;\mbox{where \emph{X} $\in$ $\mathbb{D}$}
\end{eqnarray}
are studied and the following theorem is proved
\medskip
\begin{theorem}\label{theorem-2}
There do not exist real hypersurfaces in $\mathbb{C}P^{2}$ or $\mathbb{C}H^{2}$, whose shape operator satisfies relation (\ref{Shape-Operator}).
\end{theorem}

\section{Preliminaries}
Throughout this paper all manifolds, vector fields etc are assumed to be of class $C^{\infty}$ and all manifolds are assumed to be connected. Let \emph{M} be a connected real hypersurface immersed in a non-flat complex space form $(M_{n}(c),G)$ with complex structure \emph{J} of constant holomorphic sectional curvature \emph{c}, $c\neq0$. Let \emph{N} be a locally defined unit normal vector field on \emph{M} and $\xi=-JN$. For a vector field \emph{X} tangent to \emph{M} we can write $JX=\varphi (X)+\eta(X)N$, where
$\varphi X$ and $\eta(X)N$ are the tangential and the normal component of \emph{JX} respectively. The Riemannian connection $\overline{\nabla}$ in $M_{n}(c)$ and $\nabla$ in \emph{M} are related for any vector fields \emph{X, Y} on \emph{M}
\begin{align*}
\overline{\nabla}_{Y}X=\nabla_{Y}X+g(AY,X)N,\quad\overline{\nabla}_{X}N=-AX,
\end{align*}
 where \emph{g} is the Riemannian metric on \emph{M} induced from \emph{G} of $M_{n}(c)$ and \emph{A} is the shape operator of \emph{M} in $M_{n}(c)$. \emph{M} has an almost contact metric structure $(\varphi,\xi,\eta,g)$ induced from \emph{J} on $M_{n}(c)$, where $\varphi$ is a (1,1) tensor field and $\eta$ an 1-form on \emph{M} such that
\begin{align*}
g(\varphi X,Y)=G(JX,Y),\quad\eta(X)=g(X,\xi)=G(JX,N).
\end{align*}
We thus have
\begin{eqnarray}
&&\varphi^{2}X=-X+\eta(X)\xi,\;\;\;\eta\circ\varphi=0,\;\;\; \varphi\xi=0,\;\;\;\eta(\xi)=1,\label{A1}\\
&&g(\varphi X,\varphi Y)=g(X,Y)-\eta(X)\eta(Y),\;\;g(X,\varphi Y)=-g(\varphi X,Y),\label{A2}\\
&&\nabla_{X}\xi=\varphi AX,\;\;\;(\nabla_{X}\varphi)Y=\eta(Y)AX-g(AX,Y)\xi.\label{A3}
\end{eqnarray}

Since the ambient space is of constant holomorphic sectional curvature \emph{c}, the equations of Gauss and Codazzi for any vector fields \emph{X, Y, Z} on \emph{M} are respectively given by
\begin{eqnarray}
&&R(X,Y)Z=\frac{c}{4}[g(Y,Z)X-g(X,Z)Y+g(\varphi Y ,Z)\varphi X\nonumber\\
&&-g(\varphi X,Z)\varphi Y-2g(\varphi X,Y)\varphi Z]+g(AY,Z)AX-g(AX,Z)AY,\label{A4}\\
&&(\nabla_{X}A)Y-(\nabla_{Y}A)X=\frac{c}{4}[\eta(X)\varphi Y-\eta(Y)\varphi X-2g(\varphi X,Y)\xi],\label{A5}
\end{eqnarray}
where \emph{R} denotes the Riemannian curvature tensor on \emph{M}.

Relation (\ref{A4}) implies that the \emph{structure Jacobi operator} \emph{l} is given by
\begin{eqnarray}\label{A6}
lX=\frac{c}{4}[X-\eta(X)\xi]+\alpha AX-\eta(AX)A\xi.
\end{eqnarray}

For every point \emph{P} $\in$ \emph{M}, the tangent space $T_{P}M$ can be decomposed as following
\begin{align*}
T_{P}M=span\{\xi\}\oplus \mathbb{D}
\end{align*}
where $\mathbb{D}=\ker(\eta)=\{X\;\in\;T_{P}M\mid\;\eta(X)=0\}$.

Due to the above decomposition, the vector field $A\xi$ can be written as
\begin{eqnarray}\label{A7}
A\xi=\alpha\xi+\beta U
\end{eqnarray}
where $\beta=|\varphi\nabla_{\xi}\xi|$ and $U=-\frac{1}{\beta}\varphi\nabla_{\xi}\xi\;\in\;\mathbb{D}$, provided that $\beta\neq0$.

\section{\textbf{Some Previous Results}}
Let \emph{M} be a non-Hopf hypersurface in $\mathbb{C}P^{2}$ or $\mathbb{C}H^{2}$, i.e. $M_{2}(c)$, $c\neq0$. Then the following relations hold on every three-dimensional real hypersurface in $M_{2}(c)$.
\begin{lemma}\label{lemma-1}
Let M be a real hypersurface in $M_{2}(c)$. Then the following relations hold on M
\begin{eqnarray}
&&AU=\gamma U+\delta\varphi U+\beta\xi,\;\;A\varphi U=\delta U+\mu\varphi U,\label{B1}\\
&&\nabla_{U}\xi=-\delta U+\gamma\varphi U,\;\;\nabla_{\varphi U}\xi=-\mu U+\delta\varphi U,\;\;\nabla_{\xi}\xi=\beta\varphi U,\label{B2}\\
&&\nabla_{U}U=\kappa_{1}\varphi U+\delta\xi,\;\;\nabla_{\varphi U}U=\kappa_{2}\varphi U+\mu\xi,\;\;\nabla_{\xi}U=\kappa_{3}\varphi U,\label{B3}\\
&&\nabla_{U}\varphi U=-\kappa_{1}U-\gamma\xi,\;\;\nabla_{\varphi U}\varphi U=-\kappa_{2}U-\delta\xi,\;\;\nabla_{\xi}\varphi U=-\kappa_{3}U-\beta\xi,\label{B4}
\end{eqnarray}
where $\gamma,$ $\delta,$ $\mu$, $\kappa_{1}$, $\kappa_{2}$, $\kappa_{3}$ are smooth functions on $\emph{M}$ and $\{U,\varphi U,\xi\}$ is an orthonormal basis of $\emph{M}$.
\end{lemma}
\textsc{Proof:} Since $g(AU,\xi)=g(U,A\xi)=\beta$ and $g(A\varphi U,\xi)=g(\varphi U,A\xi)=0$ we have $$\hspace{10pt}AU=\gamma U+\delta\varphi U+\beta\xi\hspace{30pt}A\varphi U=\delta U+\mu\varphi U,$$ where $\gamma,\delta,\mu$ are smooth functions.

The first relation (\ref{A3}), because of (\ref{A7}) and (\ref{B1}), for $X=U$, $X=\varphi U$ and $X=\xi$, implies (\ref{B2}).

From the well known relation $Xg(Y,Z)=g(\nabla_{X}Y,Z)+g(Y,\nabla_{X}Z)$ for \emph{X, Y, Z} $\in$ $\{\xi, U, \varphi U\}$ we obtain (\ref{B3}) and (\ref{B4}), where $\kappa_{1}, \kappa_{2}$ and $\kappa_{3}$ are smooth functions. \qed
\medskip

Owing to relation (\ref{B1}), relation (\ref{A6}) implies
\begin{equation}\label{SJO}
lU=(\frac{c}{4}+\alpha\gamma-\beta^{2})U+\alpha\delta\varphi U,\;\;\;l\varphi U=\alpha\delta U+(\alpha\mu+\frac{c}{4})\varphi U\;\;\mbox{and}\;\;l\xi=0.
\end{equation}

Because of Lemma \ref{lemma-1} the Codazzi equation (see (\ref{A5})) for $X$ $\in$ $\{U, \varphi U\}$ and $Y=\xi$ implies the following relations
\begin{eqnarray}
U\beta-\xi\gamma&=&\alpha\delta-2\delta\kappa_{3}\label{B5}\\
\xi\delta&=&\alpha\gamma+\beta\kappa_{1}+\delta^{2}+\mu\kappa_{3}+\frac{c}{4}-\gamma\mu-\gamma\kappa_{3}-\beta^{2}\label{B6}\\
U\alpha-\xi\beta&=&-3\beta\delta\label{B7}\\
\xi\mu&=&\alpha\delta+\beta\kappa_{2}-2\delta\kappa_{3}\label{B8}\\
(\varphi U)\alpha&=&\alpha\beta+\beta\kappa_{3}-3\beta\mu\label{B9}\\
(\varphi U)\beta&=&\alpha\gamma+\beta\kappa_{1}+2\delta^{2}+\frac{c}{2}-2\gamma\mu+\alpha\mu\label{B10}
\end{eqnarray}
and for $X=U$ and $Y=\varphi U$
\begin{eqnarray}
U\delta-(\varphi U)\gamma&=&\mu\kappa_{1}-\kappa_{1}\gamma-\beta\gamma-2\delta\kappa_{2}-2\beta\mu\label{B11}\\
U\mu-(\varphi U)\delta&=&\gamma\kappa_{2}+\beta\delta-\kappa_{2}\mu-2\delta\kappa_{1}\label{B12}
\end{eqnarray}

We recall the following Proposition (\cite{IR}):
\begin{proposition}\label{proposition-1}
There do not exist real hypersurfaces in $M_{2}(c)$, whose structure Jacobi opeator vanishes.
\end{proposition}

\section{\textbf{Proof of Theorem 1}}
We consider the open subset $\mathcal{W}$ of points \emph{P} $\in$ \emph{M}, such that there exists a neighborhood of every \emph{P}, where $\beta=0$ and $\mathcal{N}$ the open subset of points \emph{Q} $\in$ \emph{M},
such that there exists a neighborhood of every \emph{Q}, where $\beta\neq0$. Since, $\beta$ is a smooth function on \emph{M}, then $\mathcal{W}\cup\mathcal{N}$ is an open and dense subset of \emph{M}. In $\mathcal{W}$ $\xi$ is principal. Furthermore, we consider $\mathcal{V}$, $\Omega$ open subsets of $\mathcal{N}$
\begin{align*}
\mathcal{V}&=\{Q\;\;\in\;\;\mathcal{N}\mid\alpha=0,\;\mbox{in a neighborhood of \emph{Q}}\},\\
\Omega&=\{Q\;\;\in\;\;\mathcal{N}\mid\alpha\neq0,\;\mbox{in a neighborhood of \emph{Q}}\},
\end{align*}
where $\mathcal{V}\cup\Omega$ is open and dense in the closure of $\mathcal{N}$.
\begin{lemma}\label{lemma-2}
Let $M$ be a real hypersurface in $M_{2}(c)$, $c\neq0$, whose structure Jacobi operator satisfies relation (\ref{Structure-Jacobi-Operator}). Then the open subset $\mathcal{V}$ is empty.
\end{lemma}
\textsc{Proof:}  Relation (\ref{SJO}) on $\mathcal{V}$ becomes
\begin{align*}
lU=(\frac{c}{4}-\beta^{2})U,\;\;\;l\varphi U=\frac{c}{4}\varphi U\;\;\mbox{and}\;l\xi=0.
\end{align*}
Moreover, relation (\ref{Structure-Jacobi-Operator}) implies
\begin{eqnarray}\label{SJO1}
\nabla_{lY}X=l\nabla_{Y}X,\;\mbox{where \emph{X} $\in$ $\mathbb{D}$ and \emph{Y} $\in$ $TM$.}
\end{eqnarray}

Relation (\ref{SJO1}), because of Lemma \ref{lemma-1} and the above relation implies
\begin{eqnarray}
\kappa_{3}&=&0,\;\mbox{for $X=U$ and $Y=\xi$}\label{F1},\\
\delta=\kappa_{2}&=&0,\;\mbox{for $X=Y=\varphi U$},\\
\kappa_{1}&=&0,\;\mbox{for $X=Y=U$},\\
\mu&=&0,\;\mbox{for $X=U$ and $Y=\varphi U$}\label{F2}.
\end{eqnarray}
Relation (\ref{B6}), because of (\ref{F1})-(\ref{F2}) implies $\beta^{2}=\frac{c}{4}$. Differentiation of the last relation with respect to $\varphi U$ and taking into consideration relations (\ref{B10}) and (\ref{F1})-(\ref{F2}) leads to $c=0$, which is impossible and this completes the proof of the present Lemma.
\qed

\begin{lemma}\label{lemma-3}
Let $M$ be a real hypersurface in $M_{2}(c)$, $c\neq0$, whose structure Jacobi operator satisfies relation (\ref{Structure-Jacobi-Operator}). Then the open subset $\Omega$ is empty.
\end{lemma}
\textsc{Proof:} First of all, relation (\ref{Structure-Jacobi-Operator}) implies
\begin{eqnarray}\label{SJO2}
\nabla_{lY}X=l\nabla_{Y}X,\;\mbox{where \emph{X} $\in$ $\mathbb{D}$ and \emph{Y} $\in$ $TM$.}
\end{eqnarray}

For $X=U$ and $Y=\xi$ relation (\ref{SJO2}), because of (\ref{B3}) implies $\kappa_{3}l\varphi U=0$. We consider $\Omega_{1}$, $\Omega'_{1}$ the open subsets of $\Omega$
\begin{align*}
\Omega_{1}&=\{Q\;\;\in\;\;\Omega\mid\kappa_{3}\neq0,\;\mbox{in a neighborhood of \emph{Q}}\},\\
\Omega'_{1}&=\{Q\;\;\in\;\;\Omega\mid\kappa_{3}=0,\;\mbox{in a neighborhood of \emph{Q}}\},
\end{align*}
where $\Omega_{1}\cup\Omega'_{1}$ is open and dense in the closure of $\Omega$.

In $\Omega_{1}$ we have $l\varphi U=0$. Relation (\ref{SJO2}) for  $X=\varphi U$ and $Y=\xi$, because of (\ref{B4}) implies $\kappa_{3}lU=0$. Since $\kappa_{3}\neq0$ this results in $lU=0$ and so in $\Omega_{1}$ the structure Jacobi operator vanishes. Due to Proposition \ref{proposition-1} we obtain that the subset $\Omega_{1}$ is empty. Thus in $\Omega$ relation $\kappa_{3}=0$ holds.

Relation (\ref{SJO2}) for $X=U$ and $Y=U$ because of (\ref{B3}) and (\ref{SJO}) gives
\begin{eqnarray}
&&\delta\kappa_{1}=0,\label{C1}\\
&&\kappa_{1}(\alpha\gamma-\beta^{2}-\alpha\mu)+\alpha\delta\kappa_{2}=0,\label{C2}\\
&&\delta(\frac{c}{4}+\alpha\gamma-\beta^{2}+\alpha\mu)=0\label{C3}.
\end{eqnarray}
Because of relation (\ref{C1}), we consider $\Omega_{2}$, $\Omega'_{2}$ open subsets of $\Omega$
\begin{align*}
\Omega_{2}&=\{Q\;\;\in\;\;\Omega\mid\delta\neq0,\;\mbox{in a neighborhood of \emph{Q}}\},\\
\Omega'_{2}&=\{Q\;\;\in\;\;\Omega\mid\delta=0,\;\mbox{in a neighborhood of \emph{Q}}\},
\end{align*}
where $\Omega_{2}\cup\Omega'_{2}$ is open and dense in the closure of $\Omega$.

So in $\Omega_{2}$ we obtain $\kappa_{1}=0$ and relation (\ref{C2}) implies $\kappa_{2}=0$.
The Riemannian curvature on \emph{M} is given by the Gauss equation (see (\ref{A4})) and by the relation
\begin{eqnarray}\label{RC}
R(X,Y)Z=\nabla_{X}\nabla_{Y}Z-\nabla_{Y}\nabla_{X}Z-\nabla_{[X,Y]}Z.
\end{eqnarray}

The combination of relations (\ref{A4}) and (\ref{RC}) for $X=U$, $Y=\xi$ and $Z=\varphi U$ taking into account Lemma \ref{lemma-1}, (\ref{B5}) and  $\kappa_{1}=\kappa_{2}=\kappa_{3}=0$ implies $\delta=0$, which is a contradiction. Therefore $\Omega_{2}=\emptyset$ and in $\Omega$ relation $\delta=0$ holds.

Resuming in $\Omega$ the following relations hold
\begin{align*}
\delta=\kappa_{3}=0 \;\;\; \mbox{and relation (\ref{C2}) becomes}\;\;\;\kappa_{1}(\alpha\gamma-\beta^{2}-\alpha\mu)=0.
\end{align*}
Relation (\ref{SJO2}) owing to (\ref{B3}), (\ref{B4}) and (\ref{SJO}) yields
\begin{eqnarray}
&&\kappa_{2}(\alpha\gamma-\beta^{2}-\alpha\mu)=0,\;\mbox{for $X=Y=\varphi U$}\label{C5}\\
&&\mu(\frac{c}{4}+\alpha\mu)=0,\;\mbox{for $X=U$ and $Y=\varphi U$}.\label{C6}
\end{eqnarray}

Due to $\kappa_{1}(\alpha\gamma-\beta^{2}-\alpha\mu)=0$  we consider $\Omega_{3}$, $\Omega'_{3}$ be open subsets of $\Omega$
\begin{align*}
\Omega_{3}&=\{Q\;\;\in\;\;\Omega\mid\alpha\gamma\neq\beta^{2}+\alpha\mu,\;\mbox{in a neighborhood of \emph{Q}}\},\\
\Omega'_{3}&=\{Q\;\;\in\;\;\Omega\mid\alpha\gamma=\beta^{2}+\alpha\mu,\;\mbox{in a neighborhood of \emph{Q}}\},
\end{align*}
where $\Omega_{3}\cup\Omega'_{3}$ is open and dense in the closure of $\Omega$.

Therefore in $\Omega_{3}$ taking into account relation (\ref{C5}) we obtain $\kappa_{1}=\kappa_{2}=0$. Owing to (\ref{C6}) we consider $\Omega_{31}$, $\Omega'_{31}$ open subsets of $\Omega_{3}$
\begin{align*}
\Omega_{31}&=\{Q\;\;\in\;\;\Omega_{3}\mid\mu\neq0,\;\mbox{in a neighborhood of \emph{Q}}\},\\
\Omega'_{31}&=\{Q\;\;\in\;\;\Omega_{3}\mid\mu=0,\;\mbox{in a neighborhood of \emph{Q}}\},
\end{align*}
where $\Omega_{31}\cup\Omega'_{31}$ is open and dense in the closure of $\Omega_{3}$.

In $\Omega_{31}$ we have $\mu=-\frac{c}{4\alpha}$. The combination of (\ref{A4}) and (\ref{RC}) for $X=\varphi U$, $Y=\xi$ and $Z=U$ taking into consideration Lemma \ref{lemma-1} and (\ref{B8}) implies $c=0$, which is impossible. Thus, $\Omega_{31}$ is empty and in $\Omega_{3}$ relation $\mu=0$ holds .

In $\Omega_{3}$ the combination of (\ref{A4}) and (\ref{RC}) taking into account Lemma \ref{lemma-1} and (\ref{B12})  for $X=U$, $Y=\varphi U$ and $Z=U$ implies $c=0$, which is a contradiction. So $\Omega_{3}$ is empty.

Resuming in $\Omega$ the following relations hold
\begin{align*}
\delta=\kappa_{3}=0,\quad\alpha\gamma=\beta^{2}+\alpha\mu\;\;\mbox{and relation (\ref{C6})}.
\end{align*}
Owing to (\ref{C6}) let $\Omega_{4}$, $\Omega'_{4}$ be open subsets of $\Omega$
\begin{align*}
\Omega_{4}&=\{Q\;\;\in\;\;\Omega\mid\mu\neq0,\;\mbox{in a neighborhood of \emph{Q}}\},\\
\Omega'_{4}&=\{Q\;\;\in\;\;\Omega\mid\mu=0,\;\mbox{in a neighborhood of \emph{Q}}\},
\end{align*}
where $\Omega_{4}\cup\Omega'_{4}$ is open and dense in the closure of $\Omega$.

Then in $\Omega_{4}$ we obtain $\mu=-\frac{c}{4\alpha}$ and $\gamma=\frac{\beta^{2}}{\alpha}-\frac{c}{4\alpha}$ and relation (\ref{SJO}) leads to $lU=l\varphi U=0$. We lead to the result that the structure Jacobi operator vanishes and so because of Proposition \ref{proposition-1} we obtain $\Omega_{4}=\emptyset$.

Thus in $\Omega$ we have $\mu=0$ and $\gamma=\frac{\beta^{2}}{\alpha}$ and  relations (\ref{B6}) and (\ref{B9})-(\ref{B11}) become respectively
\begin{eqnarray}
&&\beta\kappa_{1}+\frac{c}{4}=0,\label{C7}\\
&&(\varphi U)\alpha=\alpha\beta,\label{C8}\\
&&(\varphi U)\beta=\beta^{2}+\beta\kappa_{1}+\frac{c}{2},\label{C9}\\
&&(\varphi U)\frac{\beta^{2}}{\alpha}=\frac{\beta^{2}}{\alpha}(\kappa_{1}+\beta)\label{C10}.
\end{eqnarray}
Substitution in (\ref{C10}), relations (\ref{C7})-(\ref{C9}) yields $c=0$, which is impossible. So $\Omega$ is empty and this completes the proof of the present Lemma.
\qed

From Lemmas \ref{lemma-2} and \ref{lemma-3} we lead to the following proposition
\begin{proposition}\label{proposition-2}
Let $M$ be a real hypersruface in $M_{2}(c)$, $c\neq0$, whose structure Jacobi operator satisfies relation (\ref{Structure-Jacobi-Operator}). Then M is a Hopf hypersurface.
\end{proposition}

Because of Proposition \ref{proposition-2} we have that $A\xi=\alpha\xi$ and due to Theorem 2.1 \cite{NR1}, $\alpha$ is constant. We consider a point $P$ $\in$ \emph{M} and we choose principal vector field \emph{Z} $\in$ $\ker(\eta)$ at \emph{P}, such that $AZ=\lambda Z$ and $A\varphi Z=\nu\varphi Z$. Then $\{Z, \varphi Z, \xi\}$ is a local orthonormal basis and the following relation holds (Corollary 2.3, \cite{NR1})
\begin{eqnarray}\label{C11}
\lambda\nu=\frac{\alpha}{2}(\lambda+\nu)+\frac{c}{4}.
\end{eqnarray}
The first relation of (\ref{A3}) for $X=Z$ and $X=\varphi Z$, because of $AZ=\lambda Z$ and $A\varphi Z=\nu\varphi Z$ implies
\begin{eqnarray}\label{d10}
\nabla_{Z}\xi=\lambda\varphi Z\;\;\;\mbox{and}\;\;\;\nabla_{\varphi Z}\xi=-\nu Z.
\end{eqnarray}
Relation (\ref{A6}) for \emph{X} $\in$ $\{Z,\varphi Z\}$, due to $AZ=\lambda Z$ and $A\varphi Z=\nu\varphi Z$  yields
\begin{eqnarray}\label{d11}
lZ=(\frac{c}{4}+\alpha\lambda)Z\;\;\;\mbox{and}\;\;\;l\varphi Z=(\frac{c}{4}+\alpha\nu)\varphi Z.
\end{eqnarray}
Relation (\ref{Structure-Jacobi-Operator}) taking into account(\ref{d11}) implies
\begin{align*}
(\frac{c}{4}+\alpha\nu)\nabla_{\varphi Z}Z&=l\nabla_{\varphi Z}Z,\;\;\mbox{for $X=Z$ and $Y=\varphi Z$}\\
(\frac{c}{4}+\alpha\lambda)\nabla_{Z}\varphi Z&=l\nabla_{Z}\varphi Z,\;\;\mbox{for $X=\varphi Z$ and $Y=Z$}.
\end{align*}
Taking the inner product of the latter with $\xi$, because of (\ref{d10}) we obtain respectively
\begin{eqnarray}\label{d12}
\nu(\frac{c}{4}+\alpha\nu)=0\;\;\mbox{and}\;\;\lambda(\frac{c}{4}+\alpha\lambda)=0.
\end{eqnarray}

Suppose that $\lambda$, $\nu$ are distinct at point \emph{P} $\in$ $M$. Owing to the first of (\ref{d12}), we suppose that $\nu=0$, then then the second implies $\alpha\lambda=-\frac{c}{4}$. Substituting the latter relations in (\ref{C11}) results in $c=0$, which is impossible. If we suppose that $\alpha\nu=-\frac{c}{4}$, then by following the same procedure we lead to a contradiction.

So the remaining case is that of $\lambda=\nu$ at any point \emph{P} $\in$ \emph{M}. Relation (\ref{d12}) yields that locally we have that either $\lambda=0$ or $\alpha\lambda=-\frac{c}{4}$. In both cases substitution of these relations (\ref{C11}) leads to a contradiction. With this the proof of Theorem \ref{theorem-1} is completed.

\section{\textbf{\textsc{Proof of Theorem 2}}}
We consider the open subsets $\mathcal{W}$, $\mathcal{N}$ of \emph{M}
\begin{align*}
\mathcal{W}&=\{P\;\;\in\;\;M\mid\beta=0,\;\mbox{in a neighborhood of P}\},\\
\mathcal{N}&=\{P\;\;\in\;\;M\mid\beta\neq0,\;\mbox{in a neighborhood of P}\},
\end{align*}
where $\mathcal{W}\cup\mathcal{N}$ is an open and dense subset of \emph{M}, since $\beta$ is a smooth function on \emph{M}. In $\mathcal{W}$ $\xi$ is principal. So in what follows we work in the open subset $\mathcal{N}$.

\begin{lemma}\label{lemma-4}
Let $M$ be a real hypersurface in $M_{2}(c)$, $c\neq0$, whose structure Jacobi operator satisfies relation (\ref{Shape-Operator}). Then the open subset $\mathcal{N}$ is empty.
\end{lemma}
\textsc{Proof:} Relation (\ref{Shape-Operator}) implies
\begin{eqnarray}\label{e1}
\nabla_{AY}X=A\nabla_{Y}X,\mbox{where \emph{X} $\in$ $\mathbb{D}$ and \emph{Y} $\in$ \emph{TM}.}
\end{eqnarray}

Relation (\ref{e1}) for $X=U$ and $Y=\xi$, due to relations (\ref{A7}), (\ref{B1}) and (\ref{B3}) implies
\begin{eqnarray}
&&A\nabla_{\xi}U=\nabla_{A\xi}U\Rightarrow \kappa_{3}A\varphi U=\alpha\nabla_{\xi}U+\beta\nabla_{U}U\nonumber \\
&&\Rightarrow -\kappa_{3}\delta U+(\alpha\kappa_{3}+\beta\kappa_{1}-\kappa_{3}\mu)\varphi U+\beta\delta\xi=0.\nonumber
\end{eqnarray}
so we lead to the following relations
\begin{eqnarray}\label{e2}
\delta=0\;\;\mbox{and}\;\;\kappa_{3}\mu=\alpha\kappa_{3}+\beta\kappa_{1}.
\end{eqnarray}
Relation (\ref{e1}) for $X=U$ and $Y=\varphi U$, taking into account the first of relation (\ref{e2}), (\ref{A7}), (\ref{B1}) and (\ref{B3}) gives
\begin{eqnarray}
&&A\nabla_{\varphi U}U=\nabla_{A\varphi U}U\Rightarrow \kappa_{2}A\varphi U+\mu A\xi=\mu\nabla_{\varphi U}U\nonumber\\
&&\Rightarrow \beta\mu U+\mu(\alpha-\mu)\xi=0,\nonumber\
\end{eqnarray}
from which yields
\begin{eqnarray}\label{e3}
\mu=0
\end{eqnarray}

Relation (\ref{e1}) for $X=\varphi U$ and $Y=\xi$, taking into account relations (\ref{A7}), (\ref{B1}), (\ref{B4}) (\ref{e2}) and (\ref{e3}) gives
\begin{eqnarray}
&&A\nabla_{\xi}\varphi U=\nabla_{A\xi}\varphi U\Rightarrow -\kappa_{3}AU-\beta A\xi=\alpha\nabla_{\xi}\varphi U+\beta\nabla_{U}\varphi U\nonumber \\
&&\Rightarrow (\kappa_{3}\gamma+\beta^{2})U+\beta(\kappa_{3}-\gamma)\xi=0,\nonumber\
\end{eqnarray}
which implies
\begin{eqnarray}\label{e4}
\kappa_{3}=\gamma\;\;\mbox{and}\;\;\kappa_{3}\gamma+\beta^{2}=0.
\end{eqnarray}
Substituting in the second of (\ref{e4}) the first relation of (\ref{e4}) we obtain $\beta^{2}+\gamma^{2}=0$, from which we have $\beta=0$, which is a contradiction. This completes the proof of the present Lemma
\qed

So because of Lemma \ref{lemma-4} we lead to the following proposition
\begin{proposition}\label{proposition-3}
Let $M$ be a real hypersurface in $M_{2}(c)$, $c\neq0$, whose shape operator satisfies relation (\ref{Shape-Operator}). Then M is a Hopf hypersurface.
\end{proposition}

Because of Proposition \ref{proposition-3} we have that $A\xi=\alpha\xi$ and due to Theorem 2.1 \cite{NR1}, $\alpha$ is constant. We consider a point $P$ $\in$ \emph{M} and choose principal vector field \emph{Z} $\in$ $\ker(\eta)$ at \emph{P}, such that $AZ=\lambda Z$ and $A\varphi Z=\nu\varphi Z$. Then $\{Z, \varphi Z, \xi\}$ is a local orthonormal basis and the following relation holds (Corollary 2.3 \cite{NR1})
\begin{eqnarray}\label{K1}
\lambda\nu=\frac{\alpha}{2}(\lambda+\nu)+\frac{c}{4}.
\end{eqnarray}
The first relation of (\ref{A3}) for $X=Z$ and $X=\varphi Z$ implies respectively
\begin{eqnarray}\label{K2}
\nabla_{Z}\xi=\lambda\varphi Z\;\;\mbox{and}\;\;\nabla_{\varphi Z}\xi=-\nu Z.
\end{eqnarray}
Relation (\ref{Shape-Operator}) for $X,Y$ $\in$ $\{Z,\varphi Z\}$ taking into account $AZ=\lambda Z$ and $A\varphi Z=\nu\varphi Z$ gives
\begin{align*}
A\nabla_{\varphi Z}Z&=\nu\nabla_{\varphi Z}Z,\;\;\mbox{for $X=Z$ and $Y=\varphi Z$}\\
A\nabla_{Z}\varphi Z&=\lambda\nabla_{Z}\varphi Z,\;\;\mbox{for $X=\varphi Z$ and $Y=Z$}.
\end{align*}

The inner product of the above relations with $\xi$ taking into account (\ref{K2}) implies respectively
\begin{eqnarray}\label{e7}
\nu(\nu-\alpha)=0\;\;\mbox{and}\;\;\lambda(\lambda-\alpha)=0.
\end{eqnarray}

Suppose that $\lambda$, $\nu$ are distinct at a point \emph{P}. Owing to the first of (\ref{e7}), we suppose that $\nu=0$ then the second relation yields $\lambda=\alpha$. Substituting the latter relations in (\ref{K1}) implies $c=-2\alpha^{2}$, so $c<0$ and the real hypersurface has three distinct eigenvalues. From this we conclude that the real hypersurface is of type B in $\mathbb{C}H^{2}$. Substituting the eigenvalues of this real hypersurface in relation $\nu=0$ leads to a contradiction (see \cite{Ber}). If we suppose that $\nu=\alpha$, then by following the same procedure we lead to same conclusion.

So the remaining case is that of $\lambda=\nu$ at any point $P$ $\in$ \emph{M}. Relation (\ref{e7}) implies that locally we have either $\lambda=0$ or $\lambda=\alpha$. In both cases substitution of these relations in (\ref{K1}) leads to a contradiction. Taking into account all the above, proof of Theorem \ref{theorem-2} is completed.

\section*{\textsc{Acknowledgements}} The author would like to thank Professor Philippos J. Xenos for his comments on the manuscript.

\end{document}